\documentclass[preprint]{elsarticle}

\let\today\relax
\makeatletter
\def\ps@pprintTitle{%
    \let\@oddhead\@empty
    \let\@evenhead\@empty
    \def\@oddfoot{\footnotesize\itshape
         {Submitted preprint} \hfill\today}%
    \let\@evenfoot\@oddfoot
    }
\makeatother

\usepackage{lineno,hyperref}
\usepackage{amsmath}
\usepackage{amsfonts}
\usepackage{amssymb}
\usepackage{graphicx}
\usepackage{lipsum}
\usepackage{epsfig}
\usepackage{epstopdf}
\usepackage{algorithm}
\usepackage{algpseudocode}
\usepackage{calc}
\usepackage{tikz}
\usetikzlibrary{calc,arrows,arrows.meta,decorations.markings}
\usepackage{accents}
\newcommand{\dbtilde}[1]{\accentset{\approx}{#1}}

\newtheorem{theorem}{Theorem}
\newtheorem{definition}{Definition}
\newtheorem{lemma}{Lemma}

\modulolinenumbers[5]

\begin{document}

\begin{frontmatter}

\title{Adapting free-space fast multipole method for layered media Green's function: algorithm and analysis}
\author{Min Hyung Cho}
\address{Department of Mathematical Sciences, University of Massachusetts Lowell, Lowell, MA 01854-2874}

\author{Jingfang Huang}
\address{Department of Mathematics, University of North Carolina at Chapel Hill, Chapel Hill, NC 27599-3250}
\cortext[mycorrespondingauthor]{Corresponding author}
\ead{huang@email.unc.edu}

\begin{abstract}
In this paper, we present a numerical algorithm for the accurate and efficient computation of the convolution of 
the frequency domain layered media Green's function with a given density function. Instead of compressing the 
convolution matrix directly as in the classical fast multipole method, fast direct solvers, and fast $\mathcal{H}$-matrix
algorithms, the new algorithm considers a translated form of the original matrix so that
many existing building blocks from the highly optimized free-space fast multipole method 
can be easily adapted to the Sommerfeld integral representations of the layered  
media Green's function. An asymptotic analysis is performed on the Sommerfeld integrals for large orders 
to provide an estimate of the decay rate in the 
new ``multipole" and ``local" expansions. In order to avoid the highly oscillatory integrand
in the original Sommerfeld integral representations when the source and target are close to each other, or when they
are both close to the interface in the scattered field, mathematically equivalent alternative direction integral 
representations are introduced. The convergence of the multipole
and local expansions and formulas and quadrature rules for the original and alternative direction
integral representations are numerically validated. 
\end{abstract}



\end{frontmatter}

\section{Introduction}
The efficient numerical simulation of waves in layered media in the frequency domain is a challenging task in scientific computing. One of the key numerical difficulties arises from the Sommerfeld integral representation for a general layered media Green's function \cite{chew1995waves,cai2002algorithmic}.  For example, consider a simple three-layered medium with layer interfaces at $y=0$ and $y=-d$,
and assume $k_1, k_2$, and $k_3$ are the wave numbers in each layer, respectively. 
When the source is located in the top layer $\mathbf{x}_0 = (x_0, y_0)$ with $y_0 > 0$, 
the scattered field at $\mathbf{x}_1 = (x_1, y_1)$  in the middle layer ($-d < y_1 < 0$) 
has the reflected field from the top and bottom interfaces. The reflected field 
from the top interface $g_2^t$ (formula for the bottom interface is similar) is represented by
\begin{equation}
 g_2^t(\mathbf{x}_1,\mathbf{x}_{0}) =  \int_{-\infty}^{\infty}
	e^{\sqrt{\lambda^2-k_2^{2}} y_1 }e^{i\lambda x_1} 
	e^{-\sqrt{\lambda^2-k_1^{2}} y_0 }e^{-i\lambda x_0} 
	\frac{1}{4\pi \sqrt{\lambda^2-k_2^2} } \sigma_2^t(\lambda) d\lambda \nonumber
\end{equation}
where 
{\tiny{
\begin{align}
\sigma_2^t(\lambda)=\frac{e^{d \sqrt{\lambda ^2-k_2^2}} 
(\lambda ^2+\sqrt{\lambda ^2-k_2^2} \sqrt{\lambda ^2-k_3^2}-k_2^2)}
{\sinh (d \sqrt{\lambda ^2-k_2^2}) (\lambda ^2+\sqrt{\lambda ^2-k_1^2} 
\sqrt{\lambda ^2-k_3^2}-k_2^2)+\sqrt{\lambda ^2-k_2^2} 
(\sqrt{\lambda ^2-k_1^2}+\sqrt{\lambda ^2-k_3^2}) \cosh (d \sqrt{\lambda ^2-k_2^2})} \nonumber
\end{align}
}}is referred to as the density on the top interface found by matching the continuity conditions 
at the layer interfaces~\cite{lai2014fast}. A discretization of the integral equation description of the wave field 
in layered media leads to a linear system and each entry in the coefficient matrix requires the evaluation of one or more of 
these Sommerfeld integrals \cite{chen2016accurate,chen2018accurate}, which is very expensive 
especially when the source and target are close to the same interface, e.g., when both $y_1$ and $y_0$ are close to 
the upper ($y=0$) interface and $|x_1-x_0| \gg  |y_1-y_0|$ in the function $g_2^t$. Note that
even the state-of-the-art fast direct solvers or $\cal{H}$-matrix algorithms 
\cite{greengard2009fast,hackbusch1999sparse,hackbusch2000sparse,ho2012fast} require a sufficient number of samples
of the coefficient matrix entries before the matrix can be effectively compressed for further numerical linear algebra
operations. 

This paper is aimed at {\it completely resolving} two of the many existing challenges in the
simulation of waves in layered media: (a) an optimal fast algorithm for the application of 
the layered media Green's function to a given density function, and (b) an effective 
numerical scheme to compute the layered media Green's function when the source and target 
are close to each other, or close to the media interface. Firstly, we present a general algorithm 
framework for an optimal fast solver. The new algorithm compresses a transformed version of the original matrix, and both the
expansions and translation operators are derived using mathematical analysis techniques. The transformed matrix approach made possible for the new algorithm to use many well-optimized numerical tools from existing free-space fast multipole method (FMM) with minor or no changes. This approach is different from existing FMM, fast direct solver, and fast $\mathcal{H}$-matrix algorithms which compress the coefficient matrix directly. The error analysis of the expansion is based on the asymptotic analysis of integrals which can be generalized to layered media Green's functions in three dimensions. Secondly, by introducing different integration contours for the Sommerfeld integral
representation of the layered media Green's function, the new method can also effectively 
and accurately evaluate the interactions when either the source, or target, or both are 
close to the interface between different layers. The resulting algorithm complexity is asymptotically 
optimal $\mathcal{O}(N)$ in the low frequency regime, with a prefactor close to that of the well-developed free-space 
FMM algorithms. 

To present the ideas, we restrict our attention to the two-dimensional (2-D) layered media
in this paper and organize the paper as follows. In Section~\ref{sec:algo}, we first 
summarize all the necessary building blocks for the new algorithm and present the 
pseudocode and complexity analysis. Section~\ref{sec:Green} provides several examples of these Green's functions.
In Section~\ref{sec:analysis}, we present a detailed analysis of the numerical
algorithm, including the asymptotic expansion based truncation error analysis of the 
new ``multipole" (far-field) and ``local" expansions, and the mathematically equivalent
alternative direction Sommerfeld integral representations of the original layered
media Green's function. Section~\ref{sec:numerical} presents numerical results designed to validate the algorithm analysis. We summarize our results and discuss future work in Section~\ref{sec:conclusion}.

\section{Adapting Free-space FMM for 2-D Layered Media Green's Functions}
\label{sec:algo}
We first present the algorithm framework. Detailed analysis of each algorithm component 
will be covered in Section~\ref{sec:analysis}. This section is written for readers with sufficient
knowledge of the classical FMM \cite{greengard1987fast}, which is becoming 
a standard topic in scientific computing.

\subsection{Summary of Algorithm Building Blocks}
We use notations and terminologies commonly adopted by the FMM community, and 
present the building blocks in the same order as they appear in the algorithm.

\vspace{0.1in}
{\bf \noindent Layered Media Green's Functions.}

The layered media Green's function consists of the free-space and scattered field parts. 
As the free-space Green's function is well-studied, we focus on a 
general form of the scattered field Green's function. For a target point 
$\mathbf{x} = (x,y)$ in the layer with wave number $k$ and a source point 
$\mathbf{x}_0 = (x_0,y_0)$ in the layer with wave number $k_0$, the scattered field Green's function is given by
\begin{equation}
g^s(\mathbf{x}, \mathbf{x}_0)= \int_{-\infty}^{\infty} 
    { e^{-\sqrt{\lambda^{2}-k^2}(y+d)} e^{i\lambda x} } 
    { e^{ \pm \sqrt{\lambda^2-k_0^2} y_0 }e^{-i\lambda x_0 } } 
	{ \frac{\sigma(\lambda)}{4 \pi \sqrt{\lambda^{2}-k^2}} } d\lambda
\end{equation}
where $d$ is a constant and $\sigma(\lambda)$ is  
independent of $\mathbf{x}$ and $\mathbf{x}_0$ and converges to a constant when 
$\lambda \to \pm \infty$. When the ``+" sign is used in $e^{ \pm \sqrt{\lambda^2-k_0^2} y_0}$,  
$y+d-y_0$ is assumed to be positive in order to guarantee integrability. Sample 
Sommerfeld integral representations of layered media Green's functions for acoustic waves 
and time harmonic Maxwell's equations (vector Helmholtz equations) can be found in 
\cite{lai2014fast,cho2017efficient,cui1999fast,o2014efficient} and in Section~\ref{sec:Green}.
In the following, the mathematical formulas for the ``$+$" sign are presented. 
The formulas for the ``$-$" sign can be derived in the same way and are thus omitted. 

\vspace{0.1in}
{\bf \noindent Multipole Expansion and Source-to-Multipole (S2M) Translation Operator.}
Assume there are $N$ sources located in a box centered at $\mathbf{x}_c^s = (x^s_c,y^s_c)$, 
each carries a charge $q_j$ located at $\mathbf{x}_j = (x_j,y_j)$, their contributions to the far-field location $\mathbf{x} = (x,y)$ are given by 
\begin{align}
\phi(\mathbf{x})&= \sum_{j=1}^N q_j g^s(\mathbf{x}, \mathbf{x}_j) \label{eq:direct_sum}\\
& = \int_{-\infty}^{\infty} e^{-\sqrt{\lambda^{2}-k^2}(y+d)} e^{i\lambda x}  
\left( \sum_{j=1}^N  q_j e^{ \sqrt{\lambda^2-k_0^2} y_j }e^{-i\lambda x_j } \right)  
\frac{\sigma(\lambda)}{4 \pi \sqrt{\lambda^{2}-k^2}}   d\lambda. \nonumber
\end{align}

\begin{definition}[Multipole Expansion] The far-field contributions at $\mathbf{x}$ due to charges $q_j$ located at $\mathbf{x}_j$, $j=1,2,\cdots N$ in a box centered at $\mathbf{x}_c^s$ in Eq.~(\ref{eq:direct_sum}) have the multipole expansion:
\begin{equation}
 \phi(\mathbf{x})= \sum_{p=-\infty}^{\infty} M_p \Phi_p(x,y)
  \label{eq:multipole_exp}
\end{equation}
where the multipole coefficient is given by
\begin{equation}
M_p = \sum_{j=1}^N  q_j J_p(k_0 r_j) e^{-ip \theta_j}, \label{eq:multipole_coef} 
\end{equation}
$(r_j, \theta_j)$ are the polar coordinates of $\mathbf{x}_j$ with respect to $\mathbf{x}_c^s$, $J_p$ is the Bessel function of order $p$, 
\begin{equation} \Phi_p(x,y)=
\int_{-\infty}^{\infty}  {e^{-\sqrt{\lambda^{2}-k^2}(y-y^s_c+d)} e^{i\lambda (x-x^s_c)} } 
    { \left( \frac{\lambda - \sqrt{\lambda^2-k_0^2}}{k_0} \right)^p } 
	{ \frac{\tilde{\sigma}(\lambda)}{4 \pi \sqrt{\lambda^{2}-k^2}}  } d \lambda,
	\label{eq:mpbasis}
\end{equation}
and ${ \tilde{\sigma}(\lambda) = \sigma(\lambda) e^{-\sqrt{\lambda^{2}-k^2} y^s_c+\sqrt{\lambda^2-k_0^2} y^s_c  } }$.
\end{definition}
Eq.~(\ref{eq:multipole_coef}) is derived from the Jacobi-Anger formula \cite{abramowitz1966handbook}
$$ e^{ikr\cos{\theta}} = \sum_{m=-\infty}^{\infty}i^m J_m(kr)e^{im\theta}$$
and the multipole coefficient $M_p$ is the same as the multipole coefficient for the free-space Helmholtz kernel \cite{rokhlin1990rapid}.  We refer to forming the multipole expansion as the Source-to-Multipole (S2M) operator. This definition simply states that by changing the basis (matrix transformation) to $\Phi_p$, existing free-space S2M operator for the Helmholtz equation with the same wave number $k_0$ can be 
used directly to derive the coefficients $M_p$ of the compressed representation in 
Eq.~(\ref{eq:multipole_exp}). Note that the asymptotic properties of 
$\tilde{\sigma}(\lambda) = \sigma(\lambda) e^{-\sqrt{\lambda^{2}-k^2} y^s_c+\sqrt{\lambda^2-k_0^2} y^s_c  }$
remain the same as the original $\sigma(\lambda)$ as $\lambda \to \pm \infty$.

{\noindent \bf Remark:} By introducing the change of variable $z=\frac{\lambda - \sqrt{\lambda^2-k_0^2}}{k_0}$ in Eq.~(\ref{eq:mpbasis}),  the expansion $\sum_{p=-\infty}^{\infty} M_p z^p $ is the Laurent expansion of the function 
\[ 
 \sum_{j=1}^N  q_j e^{ \sqrt{\lambda^2-k_0^2} (y_j-y^s_c) }e^{-i\lambda (x_j-x^s_c) } 
\]
and $M_p$ is the expansion coefficient independent of $\lambda$ (or $z$).

At the low-frequency regime when both $k_0$ and $k$ are small, the number of terms required 
in the truncated layered media multipole expansion for a prescribed error tolerance is approximately the same as 
that in the free-space Laplace FMM. We leave the truncation error 
analysis of the multipole (and local) expansions to Section \ref{sec:analysis}.  

\vspace{0.1in}
{\bf \noindent Multipole-to-Multipole (M2M) Translation Operator.}
In FMM, multipole expansion of the parent is constructed by translating that of the children. Translating the center of a multipole expansion from child box to its parent box is referred to as the Multipole-to-Multipole (M2M) translation operator. As the multipole coefficients for the layered media in Eq.~(\ref{eq:multipole_coef}) are the same as those for the free-space Green's function for the Helmholtz equation with wave number $k_0$ for both the parent and child boxes, we therefore have the following lemma.

\begin{lemma}[M2M]
The M2M translation operator for the layered media Green's function in Eq.~(\ref{eq:multipole_coef}) is 
the same as the M2M operator for the free-space Green's function of the Helmholtz equation with 
wave number $k_0$. The parent's multipole coefficients $\widetilde{M}_p$ are given by 
\begin{equation}
\widetilde{M}_p = \sum_{q=-\infty}^{\infty} M_{p-q} J_q(k_0 r_{12}) e^{i q \theta_{12}}
\end{equation} 
where $M_{p-q}$ are the child's multipole coefficients and  $(r_{12}, \theta_{12})$ are the
polar coordinates of the child's center with respect to the parent's center.
\end{lemma}
Therefore, the free-space M2M translation operator can be used without any change to obtain the multipole
expansions for all the boxes on the hierarchical tree structure for the {\it layered media Green's function}.

\vspace{0.1in}
{\bf \noindent Local Expansion and Multipole-to-Local (M2L) Translation Operator.}
Notice that the potential field $\phi(\mathbf{x})$ satisfies the Helmholtz equation with wave number
$k$, we therefore use the same Bessel function based expansion as that for the free-space Green's 
function with wave number $k$ to compress the received far-field contributions into a local 
expansion of the target box centered at $\mathbf{x}_c^t = (x_c^t,y_c^t)$.
\begin{definition}[Local Expansion]
The potential function $\phi$ due to the far-field source contributions can be
	compressed into a local expansion
\begin{equation}
\phi(\mathbf{x}) = \sum_{p=-\infty}^{\infty} L_{p} J_{p} (k r)e^{i p \theta}
\label{eq:local}
\end{equation}
where $(r,\theta)$ are the polar coordinates of $\mathbf{x}$ with respect to the target box center
$\mathbf{x}_c^t$, $L_p$ is called the local expansion coefficient. Evaluating the local expansion at a target point is referred to as the Local-to-Target (L2T) translation operator.
\end{definition}

Similar to the free-space FMM, the compressed far-field multipole expansion of the source
box centered at $\mathbf{x}_c^s$ given by Eq.~(\ref{eq:multipole_exp}) can be translated into a local 
expansion in Eq.~(\ref{eq:local}) of the target box centered at $\mathbf{x}_c^t$, by 
plugging the plane wave formula 
\begin{equation}
e^{i (\lambda (x-x_c^t) + \sqrt{k^2-\lambda^2} (y-y_c^t) )} = \sum_{m=-\infty}^{\infty} i^m J_m(k r)e^{i m \theta} 
  \left( \frac{\lambda-i \sqrt{k^2-\lambda^2}}{k} \right)^m
\end{equation}
in the basis function $\Phi_p$ in Eq.~(\ref{eq:mpbasis}),
where $(r,\theta)$ are the polar coordinates of $\mathbf{x}$ with respect to the target box center at 
$\mathbf{x}_c^t$. Translating the multipole expansion to a local expansion is referred to as the Multipole-to-Local (M2L) translation operator and we have the following lemma.
\begin{lemma}[M2L]
The local expansion coefficients of a target box centered at $\mathbf{x}_c^t$, 
due to the contributions from particles in a source box centered at $\mathbf{x}_c^s$ described by its
multipole expansion in Eq.~(\ref{eq:multipole_exp}), can be computed using the M2L translation matrix $\mathrm{A} = \{A_{p,q}\}$ using
\begin{equation}
 L_p = \sum_{q=-\infty}^{\infty} A_{p,q} M_q
\label{eq:M2L}
\end{equation}
where
\begin{align}
 A_{p,q} = \int_{-\infty}^{\infty} i^p & e^{-\sqrt{\lambda^{2}-k^2}(y_c^t-y^s_c+d)} e^{i\lambda (x_c^t-x^s_c)} 
    \\
 &   \times \left(\frac{\lambda+\sqrt{ \lambda^2- k^2}}{k}\right)^p  \left(\frac{\lambda - \sqrt{\lambda^2-k_0^2}}{k_0}\right)^q  
 \frac{\tilde{\sigma}(\lambda)}{4 \pi \sqrt{\lambda^{2}-k^2}}   d\lambda. \nonumber
\end{align}
\end{lemma}
Note that when $|y_c^t-y^s_c+d| \ll |x_c^t-x^s_c|$, the evaluation of $A_{p,q}$ can be very costly due
to the highly oscillatory term $e^{i\lambda (x_c^t-x^s_c)}$. This difficulty can be resolved by
using the {\it alternative direction} integral representation of the translation operator that will be discussed in 
Sec.~\ref{sec:analysis}. The translation matrix $\mathrm{A}$ can be precomputed for optimal efficiency
when the geometries of the layered media and embedded objects are fixed, or computed on the fly 
using the Gauss and Laguerre quadrature rules for the {\it alternative direction} integral representations.

\vspace{0.1in}
{\bf \noindent Local-to-Local (L2L) Translation Operator.} The local expansion of the parent can be translated to its children 
using the Local-to-Local (L2L) translation operator given by the following lemma \cite{rokhlin1990rapid,crutchfield2006remarks}.
\begin{lemma}[L2L]
As the basis for the local expansion of the layered media Green's function is the same as that
for the free-space Green's function of the Helmholtz equation with wave number $k$, the L2L 
translation operator for the layered media Green's function in Eq.~(\ref{eq:local}) is 
therefore the same as the L2L operator of the free-space FMM for the Helmholtz equation 
with wave number $k$. The child's local coefficients $\widetilde{L}_p$ are given by 
\begin{equation}
\widetilde{L}_p = \sum_{q=-\infty}^{\infty} L_{p-q} J_q(k r_{12}) e^{ -i q (\theta_{12} -\pi) }
\end{equation} 
where $L_{p-q}$ are the parent's local expansion coefficients and  $(r_{12}, \theta_{12})$ are the
polar coordinates of the parent's box center with respect to the child's box center.
\end{lemma}
This lemma states that the free-space L2L translation operator can be used without change 
to derive the local expansion of the child box from its parent's for the layered media case.


\vspace{0.1in}
{\bf \noindent Direct Source-to-Target (S2T) Interactions.}
For neighboring boxes, the interaction of the source and target can be handled in two ways: 
(1) by evaluating the Sommerfeld integral directly, or (2) for the scattered field contribution 
from sources in a source box, when its multipole expansion is also valid in its neighboring
target box, the source box's multipole expansion can be translated and merged into the target 
box's local expansion using the same M2L translation operator as we discussed previously for 
the interaction list boxes, which will be evaluated later using the very efficient L2T 
operator for the leaf boxes. Note that when this happens, the FMM tree structure should be modified 
accordingly to further accelerate the computation.

One numerical difficulty of evaluating the layered media Green's function for direct S2T interaction or entries of the translation matrix $\mathrm{A}$ in M2L comes from the oscillatory term $e^{i \lambda x }$ when the other exponential terms in the integrand decay slowly. 
For example, in the free-space Green's function, 
\begin{align}
g(\mathbf{x},\mathbf{x}_{0})=\int_{-\infty}^{\infty}
  e^{-\sqrt{\lambda^{2}-{k}^{2}} y} e^{i\lambda x} 
 \frac{1}{4\pi \sqrt{\lambda^{2}-{k}^{2}} }  d\lambda 
	\label{eq:freespace_sommerfeld}
\end{align}
where the source point is located at the origin $(x_0,y_0)=(0,0)$ and the target point is located at the first quadrant ($x>0$ and $y>0$).  In the new version of FMM
\cite{crutchfield2006remarks,greengard1998accelerating,greengard1997new}, it is referred to as the ``north" plane wave expansion because this formula is valid when $y>0$. However, when $y$ is very close to the line $y=0$ and $x\gg 0$, a huge number of quadrature points has to 
be used to resolve the oscillatory term $e^{i\lambda x}$ due to the very slow decay of  $e^{-\sqrt{\lambda^{2}-{k}^{2}} y}$. Similar problems arise when evaluating the direct interactions of the source and target points that are close to the interface of the layered media and the M2L translation matrix $\mathrm{A}$ in Eq.~(\ref{eq:M2L}).

To understand the origin of this problem, we divide the Sommerfeld integral representation into the
propagating ($|\lambda|<k$) and evanescent ($|\lambda|<k$) parts, and after the change of variables as in \cite{greengard1998accelerating}, Eq. (\ref{eq:freespace_sommerfeld}) can be rewritten as
\begin{displaymath}
g(\mathbf{x},\mathbf{x}_{0})=  	\frac{1}{4\pi} \int_{0}^{\pi}
	e^{i k (y \sin \theta - x \cos \theta)}  d\theta + 
	\frac{1}{2\pi i}\int_{0}^{\infty}  e^{- t y} \cos(\sqrt{t^2+k^2} x)   dt.
\end{displaymath}
The first and second integrals are called the propagating part and evanescent part, respectively. For the propagating part, the number of required Gauss 
quadrature points to evaluate the integral only depends on $k$ and $r=\sqrt{x^2+y^2}$, 
which is normal. However, when $y\ll x$, the integrand in the evanescent part 
$ e^{- t y} \cos(\sqrt{t^2+k^2} x)$ requires many Laguerre quadrature points to resolve the oscillatory term $\cos(\sqrt{t^2+k^2} x)$ because $e^{- ty}$ decays
slowly. This problem can be resolved by using the equivalent ``east" plane wave representation for 
$x>0$  
\begin{displaymath}
g(\mathbf{x},\mathbf{x}_{0})=  	\frac{1}{4\pi}\int_{0}^{\pi}
	 e^{i k (-y \cos \theta + x \sin \theta)}  d\theta + 
	\frac{1}{2\pi i} \int_{0}^{\infty}  e^{- t x} \cos(\sqrt{t^2+k^2} y)   dt,
\end{displaymath}
which can be derived from the ``north" plane wave representation using contour integration as will be discussed in Section \ref{sec:analysis}.
Since $y \ll x$, the evanescent part of the alternative direction integral can be evaluated using a small number of Laguerre quadrature points. Similar 
representations can also be derived for the ``south" and ``west" directions for general layered media Green's functions, either using the method of images, contour integration, or applying the integral transforms directly in the alternative direction. 

\subsection{Algorithm Pseudocode}
We present the algorithm pseudocode in Algorithm 1. Compared with the classical FMM for free-space kernels,
the adapted FMM for layered media Green's function only differs in the M2L and S2T subroutines, and in 
the number of expansion terms when the wave numbers are different. All other subroutines and functions 
from free-space FMM can be adopted with minor or no changes by the layered media FMM.
\begin{algorithm}[htbp]
\begin{algorithmic}
\caption{Adapting Free-space FMM for Layered Media Green's Functions}

\State  {\begin{center}{\bf Step 1: Initialization} \end{center}}
\State {{\bf Generate} an adaptive hierarchical tree and precompute necessary tables.}

\State{\begin{center}{\bf Step 2: Upward Pass}\end{center}}
\For{$l = L, \cdots, 0$}
	\For{ all boxes $j$ on level $l$}
		\If {$j$ is a leaf node}
			\State{{\bf compute} S2M using {\it free-space} S2M operator.}
		\Else
			\State{{\bf compute} M2M using {\it free-space} M2M operator.}
		\EndIf
	\EndFor
\EndFor

\State{\begin{center}{\bf Step 3: Downward Pass}\end{center}}
\For{$l=1, \cdots, L$}
	\For{all boxes $j$ on level $l$}
		\State {{\bf shift} local expansion of $j$'s parent to $j$ using {\it free-space} L2L operator.}
		\State{{\bf collect} interaction list contribution using M2L operator in Eq.~(\ref{eq:M2L}).} 
		\State{{\bf collect} valid neighbor box multipole expansion using M2L operator.} 
	\EndFor	
\EndFor

\State{\begin{center}{\bf Step 4: Evaluate Local Expansions and Direct Interactions}\end{center}}
\For{each leaf node (childless box)}
	\State{{\bf evaluate} local expansion (L2T) at each particle location.}
	\State{{\bf collect} un-evaluated source target interaction (S2T) from neighbor boxes (including self) using alternative direction Sommerfeld integrals.}
\EndFor
\end{algorithmic}
\end{algorithm}

\subsection{Algorithm Complexity}
We compare the algorithm complexity of the layered media FMM (LM-FMM) with that of the free-space FMM (FS-FMM).

In the upward pass, only the free-space S2M and M2M translation operators are used. Thus, if the number of expansion terms is the same as that of the free-space case, the LM-FMM has the same number of operations as the FS-FMM in the upward pass. In the downward pass, the L2L operator of the LM-FMM has the same complexity as that of the FS-FMM. However, the M2L operator requires more operations when the M2L translation matrix is computed on the fly, as at least one integral has to be evaluated to find each entry of the translation matrix. Note that these integrals can be evaluated efficiently when the alternative direction integral representations are used (further discussed in Sec.~\ref{sec:analysis}). Moreover, the translation matrix can often be re-used by many boxes in the same level in the hierarchical tree structure. On the other hand, when the translation matrix $\mathrm{A}$ is precomputed, the algorithm complexity in the downward pass is about the 
same as that in the classical FS-FMM. In Step 4, the evaluation of the local expansion in the LM-FMM has the 
same complexity as that in the FS-FMM. For the direct source to target (S2T) interactions, more operations are needed than FS-FMM because the alternative direction Sommerfeld 
integrals have to be evaluated. In summary, when the local direct interaction operations are not counted (which can be very efficient on parallel computers), the translation matrix $\mathrm{A}$ is precomputed, and the number of expansion terms is the same as that of the FS-FMM, the total number of operations of the LM-FMM is about the same as that of the FS-FMM. In Sec. \ref{sec:analysis}, we show that for many settings, the number of expansion terms of the LM-FMM is in the same order as that of the FS-FMM, for example, when the wave numbers of different layers are all in the low-frequency regime.

\section{Examples of 2-D Layered Media Green's Functions}
\label{sec:Green}
We consider the potential at a target point $\mathbf{x}=(x,y)$ due to a source charge with density $q_0$ located at 
$\mathbf{x}_0=(x_0,y_0)$. The source and target may locate in the same or different layers with the $j^{th}$ interface 
located at $y=y_j$. The wave numbers are $k_0$ for the source layer and $k$ for the target layer,
respectively. 

To derive the layered media Green's function, a Fourier transform is usually performed in the $x$-direction,
reducing the 2-D Helmholtz equation to an ODE system, which can be solved analytically with some unknown
density functions in the Sommerfeld integral representations. The density functions are obtained by 
solving a linear system of algebraic equations to match the interface conditions. In this paper, we 
consider the Green's function in a rather general form
\begin{equation}
G(\mathbf{x}, \mathbf{x}_0)= \int_{-\infty}^{\infty} 
e^{-\sqrt{\lambda^{2}-k^2}(y+d)} e^{i\lambda x}  
e^{ \pm \sqrt{\lambda^2-k_0^2} y_0 }e^{-i\lambda x_0 }  
\frac{\sigma(\lambda)}{4 \pi \sqrt{\lambda^{2}-k^2}}  d\lambda
   \label{eq:Green}
\end{equation}
where $d$ is a constant and $\sigma(\lambda)$ converges to a constant 
when $\lambda \to \pm \infty$. Many 2-D layered media Green's functions are 
in this particular form. In the following discussions, we refer to the first two exponential terms $e^{-\sqrt{\lambda^{2}-k^2}(y+d)} e^{i\lambda x}$ as the target term, the third and fourth exponential terms $e^{ \pm \sqrt{\lambda^2-k_0^2} y_0 }e^{-i\lambda x_0 }$
as the source term, and $\sigma(\lambda)$ as the image term (for reasons which will be explained later).

In this section, we present a few examples of the layered media Green's functions, all are 
in the form of Eq.~(\ref{eq:Green}). 
 
\vspace{0.1in}
{\noindent \bf Example 1: Free-space Green's Function.} 
The first example is the free-space Green's function for the Helmholtz equation with
wave number $k$. For a source point $\mathbf{x}_0 = (x_0, y_0)$ and a target point $\mathbf{x}=(x,y)$ with $y-y_0>0$, the free-space Green's function is given by the Sommerfeld integral of the form
\begin{equation}
g(\mathbf{x},\mathbf{x}_{0})=\int_{-\infty}^{\infty}
e^{-\sqrt{\lambda^{2}-{k}^{2}} y} e^{i\lambda x} 
e^{\sqrt{\lambda^{2}-{k}^{2}} y_0} e^{- i\lambda x_0} 
\frac{1}{4\pi \sqrt{\lambda^{2}-{k}^{2}} } 
d\lambda.
\label{eq:spectralfree}%
\end{equation}

\vspace{0.1in}
{\noindent \bf Example 2: Two-layered Medium with Zero Dirichlet Interface Condition.} 
The second example is the half-space problem with Dirichlet condition (total field is $0$) at the layer interface 
located at $y=0$. We assume the source is located at $\mathbf{x}_0 = (x_0,y_0)$ on the upper half plane ($y_0>0$).
Using the method of images, the scattered field for a target point located at $\mathbf{x} = (x,y)$ on the upper
half plane ($y>0$) can be represented as
\begin{equation}
g^s(\mathbf{x},\mathbf{x}_{0})=\int_{-\infty}^{\infty}
e^{-\sqrt{\lambda^{2}-{k}^{2}} y} e^{i\lambda x} 
e^{-\sqrt{\lambda^{2}-{k}^{2}} y_0} e^{- i\lambda x_0} 
\frac{-1}{4\pi \sqrt{\lambda^{2}-{k}^{2}} } 
d\lambda.
\label{eq:dirichlet}%
\end{equation}
To satisfy the Dirichlet boundary condition, an image charge is added at the location $(x_0,-y_0)$ 
and $-1$ is the charge of the image source. 
Interested readers are referred to \cite{lai2018new} for more details.

\vspace{0.1in}
{\noindent \bf Example 3: Two-layered Medium with Impedance Interface Condition.} 
The third example is the half-space Green's function with the impedance boundary condition
\begin{equation}
\frac{\partial u}{\partial y}-i\alpha u=0
\label{impedBC}
\end{equation}
 at the layer interface $y=0$. The scattered field at a target point $\mathbf{x} =(x,y) $ due to the source at $\mathbf{x}_0=(x_0,y_0)$, where both points are located at the upper half plane with wave number $k$, is given by 
\begin{equation}
	g^s(\mathbf{x},\mathbf{x}_{0})=\int_{-\infty}^{\infty}
e^{-\sqrt{\lambda^2-k^{2}} y }e^{i\lambda x} 
e^{-\sqrt{\lambda^2-k^{2}} y_0 }e^{-i\lambda x_0} 
 \frac{1}{4\pi \sqrt{\lambda^2-k^2} } 
  \frac{\sqrt{\lambda^2-k^2}+i k \alpha}{\sqrt{\lambda^2-k^2}-i k \alpha}  d\lambda.
\label{eq:halfspace}
\end{equation}
In \cite{o2014efficient,cho2018heterogeneous}, it was shown that the term $\frac{\sqrt{\lambda^2-k^2}+i k \alpha}{\sqrt{\lambda^2-k^2}-i k \alpha}$
can be derived using the method of complex images, we therefore refer to  $\sigma(\lambda)=\frac{\sqrt{\lambda^2-k^2}+i k \alpha}{\sqrt{\lambda^2-k^2}-i k \alpha}$ term as the image term. This Green's function was also discussed in \cite{lai2018new}. Note that when 
$\lambda \to \pm \infty$, $\sigma \to 1$.

\vspace{0.1in}
{\noindent \bf Example 4: Three-layered Medium with Transmission Condition.} 
The last example we consider in this paper is the Green's function for a three-layered medium with the transmission condition in \cite{lai2014fast},
where the layer interfaces are located at $y=0$ and $y=-d$. Let a source point be located in the top layer at $\mathbf{x}_0 = (x_0, y_0)$ and a target point $\mathbf{x} = (x,y)$. In the first layer ($y > 0$),  the scattered field $g_1^s$ is the field reflected from bottom layers. In the second layer ($-d <y<0$), the scattered field consists of the contribution from top and bottom interfaces $g_2^s = g_2^t+g_2^b$. In the third layer ($y<-d$), the scattered field $g_3^s$ is the transmitted field from the source in the first layer. By matching the continuity of the field, the scattered field in each layer can be represented as
\begin{align}
&g_1^s(\mathbf{x},\mathbf{x}_{0}) =  \int_{-\infty}^{\infty}
e^{-\sqrt{\lambda^2-k_1^{2}} y }e^{i\lambda x} 
	e^{-\sqrt{\lambda^2-k_1^{2}} y_0 }e^{-i\lambda x_0} 
	\frac{1}{4\pi \sqrt{\lambda^2-k_1^2} } \sigma_1(\lambda) d\lambda, \nonumber\\
&g_2^t(\mathbf{x},\mathbf{x}_{0}) =  \int_{-\infty}^{\infty}
 e^{\sqrt{\lambda^2-k_2^{2}} y }e^{i\lambda x} 
e^{-\sqrt{\lambda^2-k_1^{2}} y_0 }e^{-i\lambda x_0} 
\frac{1}{4\pi \sqrt{\lambda^2-k_2^2} } \sigma_2^t(\lambda) d\lambda,  \nonumber\\
&g_2^b(\mathbf{x},\mathbf{x}_{0}) = \int_{-\infty}^{\infty}
 e^{-\sqrt{\lambda^2-k_2^{2}} (y+2d) }e^{i\lambda x} 
	e^{-\sqrt{\lambda^2-k_1^{2}} y_0 }e^{-i\lambda x_0} 
	\frac{1}{4\pi \sqrt{\lambda^2-k_2^2} } \sigma_2^b(\lambda) d\lambda,  \nonumber\\
&g_3^s(\mathbf{x},\mathbf{x}_{0}) =  \int_{-\infty}^{\infty}
e^{-\sqrt{\lambda^2-k_3^{2}} (y+2d) }e^{i\lambda x} 
	e^{-\sqrt{\lambda^2-k_1^{2}} y_0 }e^{-i\lambda x_0} 
	\frac{1}{4\pi \sqrt{\lambda^2-k_3^2} } \sigma_3(\lambda) d\lambda  \nonumber
\end{align}
where $(\sigma_1(\lambda), \sigma_2^t(\lambda), \sigma_2^b(\lambda), \sigma_3(\lambda))^T$ are the solutions of the linear system
\begin{equation}
\left(
\begin{array}{rrrr}
 -1 & \frac{\sqrt{\lambda ^2-k_1^2}}{\sqrt{\lambda ^2-k_2^2}} & \frac{e^{-d \sqrt{\lambda ^2-k_2^2}} \sqrt{\lambda ^2-k_1^2}}{\sqrt{\lambda ^2-k_2^2}} & 0 \\
 0 & e^{-d \sqrt{\lambda ^2-k_2^2}} & 1 & -\frac{\sqrt{\lambda ^2-k_2^2}}{\sqrt{\lambda ^2-k_3^2}} \\
 1 & 1 & -e^{-d \sqrt{\lambda ^2-k_2^2}} & 0 \\
 0 & e^{-d \sqrt{\lambda ^2-k_2^2}} & -1 & -1 \\
\end{array}
	\right) \left( 
	\begin{array}{c}
	 \sigma_1(\lambda) \\  \sigma_2^t(\lambda)\\  \sigma_2^b(\lambda)\\ \sigma_3(\lambda) 
	\end{array}
	\right) =  \left(
	\begin{array}{c}
	 1 \\  0 \\  1 \\  0
	\end{array}
	\right),
\end{equation}
and the solutions $( \sigma_1(\lambda),\sigma_2^t(\lambda),\sigma_2^b(\lambda),\sigma_3(\lambda))^T$ are explicitly given by
{\tiny{
\[
\left( \begin{array}{l}
	\frac{\sinh \left(d \sqrt{\lambda ^2-k_2^2}\right) \left(-\lambda ^2+\sqrt{\lambda ^2-k_1^2} \sqrt{\lambda ^2-k_3^2}+k_2^2\right)+\sqrt{\lambda ^2-k_2^2} \left(\sqrt{\lambda ^2-k_1^2}-\sqrt{\lambda ^2-k_3^2}\right) \cosh \left(d \sqrt{\lambda ^2-k_2^2}\right)}{\sinh \left(d \sqrt{\lambda ^2-k_2^2}\right) \left(\lambda ^2+\sqrt{\lambda ^2-k_1^2} \sqrt{\lambda ^2-k_3^2}-k_2^2\right)+\sqrt{\lambda ^2-k_2^2} \left(\sqrt{\lambda ^2-k_1^2}+\sqrt{\lambda ^2-k_3^2}\right) \cosh \left(d \sqrt{\lambda ^2-k_2^2}\right)} \\  
		\frac{e^{d \sqrt{\lambda ^2-k_2^2}} \left(\lambda ^2+\sqrt{\lambda ^2-k_2^2} \sqrt{\lambda ^2-k_3^2}-k_2^2\right)}{\sinh \left(d \sqrt{\lambda ^2-k_2^2}\right) \left(\lambda ^2+\sqrt{\lambda ^2-k_1^2} \sqrt{\lambda ^2-k_3^2}-k_2^2\right)+\sqrt{\lambda ^2-k_2^2} \left(\sqrt{\lambda ^2-k_1^2}+\sqrt{\lambda ^2-k_3^2}\right) \cosh \left(d \sqrt{\lambda ^2-k_2^2}\right)} \\  
		-\frac{e^{d \sqrt{\lambda ^2-k_2^2}} \left(-\lambda ^2+\sqrt{\lambda ^2-k_2^2} \sqrt{\lambda ^2-k_3^2}+k_2^2\right)}{\sinh \left(d \sqrt{\lambda ^2-k_2^2}\right) \left(\lambda ^2+\sqrt{\lambda ^2-k_1^2} \sqrt{\lambda ^2-k_3^2}-k_2^2\right)+\sqrt{\lambda ^2-k_2^2} \left(\sqrt{\lambda ^2-k_1^2}+\sqrt{\lambda ^2-k_3^2}\right) \cosh \left(d \sqrt{\lambda ^2-k_2^2}\right)} \\  
		\frac{2 \sqrt{\lambda ^2-k_2^2} \sqrt{\lambda ^2-k_3^2} e^{d \sqrt{\lambda ^2-k_3^2}}}{\sinh \left(d \sqrt{\lambda ^2-k_2^2}\right) \left(\lambda ^2+\sqrt{\lambda ^2-k_1^2} \sqrt{\lambda ^2-k_3^2}-k_2^2\right)+\sqrt{\lambda ^2-k_2^2} \left(\sqrt{\lambda ^2-k_1^2}+\sqrt{\lambda ^2-k_3^2}\right) \cosh \left(d \sqrt{\lambda ^2-k_2^2}\right)}
\end{array}
\right).
\]
}}
Note that when $\lambda \to \pm \infty$, all the $\sigma$ functions converge to constants. 

The form of the Green's function in Eq.~(\ref{eq:Green}) is not surprising. Clearly the first target term satisfies the Helmholtz equation
at the target layer with wave number $k$ and the source term satisfies the Helmholtz equation at the 
source layer with wave number $k_0$. The third term is independent of the variables $\mathbf{x}_0$ and 
$\mathbf{x}$, and we collect all the exponential growth or decay terms in the constant $d$ in the target term so the 
image term converges to a constant when $\lambda \to \pm \infty$.

\section{Analysis of Layered Media Fast Multipole Method}
\label{sec:analysis}
We present a detailed analysis of the algorithm for layered media Green's function in this section. 
We focus on the following two topics: (a) the number of terms in the multipole and local expansions and 
truncation errors, and (b) evaluation of the local direct interactions and M2L translation operators using 
the mathematically equivalent alternative direction Sommerfeld integral representation.

\subsection{Truncating the Multipole and Local Expansions}
\label{sec:truncate}
We study the decay rates of the terms in the multipole and local expansions by considering the setting of a single source with a unit charge located at $\mathbf{x}_0 = (x_0,y_0)$. The multipole expansion describes the 
potential at far-field locations as a function of $\mathbf{x} = (x,y)$ due to a charge in the source box, namely,
$$ \phi(\mathbf{x}) = \sum_{p=-\infty}^{\infty} M_p \Phi_p(x,y),  \quad M_p=J_p(k_0 r) e^{-ip \theta} $$ 
where $(r, \theta)$ are the polar coordinates of the point $(x_0,y_0)$ with respect to the source box center
$\mathbf{x}_c^s = (x_c^s,y_c^s)$, the basis function $\Phi_p(\mathbf{x})$ is given by
\begin{equation} 
\Phi_p(\mathbf{x})=  
\int_{-\infty}^{\infty}  
    e^{-\sqrt{\lambda^{2}-k^2}(y-y^s_c+d)} e^{i\lambda (x-x^s_c)}  
   \left( \frac{\lambda - \sqrt{\lambda^2-k_0^2}}{k_0} \right)^p  
   \frac{\tilde{\sigma} (\lambda)}{4 \pi \sqrt{\lambda^{2}-k^2}}   d \lambda
\end{equation}
and we define the ``modified" distance between $\mathbf{x}$ and $\mathbf{x}_c^s$ of the multipole
expansion as $\rho=\sqrt{(y-y^s_c+d)^2+(x-x^s_c)^2}$.
The multipole expansion can be considered as a compressed representation of the source box's contribution 
to be sent to the far-field locations. Similarly, the local expansion 
associated with each target box compresses the received far-field contributions and describes the potential 
as a function of $(x-x_c^t,y-y_c^t)$ given by
$$ \phi(\mathbf{x}) = \sum_{p=-\infty}^{\infty} L_p  { J_p(k \tilde{r}) e^{i p \tilde{\theta}}} 
  =\sum_{p=-\infty}^{\infty} { J_p(k \tilde{r}) e^{i p \tilde{\theta}}} \Psi_p(x_0,y_0)$$ 
where $(\tilde{r}, \tilde{\theta})$ are the polar coordinates of the point $(x,y)$ with respect to the target box center
$\mathbf{x}_c^t = (x_c^t,y_c^t)$ and the basis function $\Psi_p(\mathbf{x}_0)$ is given by
\begin{equation} 
\Psi_p(\mathbf{x}_0)= 
\int_{-\infty}^{\infty}  
        \left( \frac{\lambda - i \sqrt{k^2-\lambda^2}}{k} \right)^p  
	e^{-\sqrt{\lambda^{2}-k_0^2}(y^t_c+d-y_0)} e^{i\lambda (x^t_c-x_0)}  
	\frac{\dbtilde{\sigma} (\lambda)}{4 \pi \sqrt{\lambda^{2}-k^2}}  d \lambda
\label{eq:localbasis}
\end{equation}
where
$ {\dbtilde{\sigma}(\lambda) = \sigma(\lambda) e^{-\sqrt{\lambda^{2}-k^2} y^t_c+\sqrt{\lambda^2-k_0^2} y^t_c  }. }$
We define the ``modified" distance between $\mathbf{x}_0$ and $\mathbf{x}_c^t$ of the local
expansion as $\tilde{\rho}=\sqrt{(y^t_c-y_0+d)^2+(x^t_c-x_0)^2}$.

Similar to the FS-FMM where the truncation errors of the multipole and local expansions are 
determined by the decay rate of $|J_p(kr) H_p(k \rho)|$, we study the decay rate of the term 
$|J_p(k_0 r) \Phi_p(x,y)|$ for the multipole expansion, and $|J_p(k r) \Psi_p(x_0,y_0)|$ for 
the local expansion for different physical parameter settings. When all other variables are 
fixed, $J_p \to 0$ and $H_p,\Phi_p,\Psi_p \to \infty$ as $p \to \infty$. Therefore, it is necessary to understand the asymptotic behavior of these functions for large $p$ values. 

\vspace{0.1in}
{\bf \noindent Asymptotic Forms of Bessel Functions for Large Order.} The asymptotic expansion for large order Bessel
functions is a well-studied topic. We cite the following well-known results from \cite{abramowitz1966handbook}, which
are valid for fixed $z$ when $\nu \to \infty$.
\begin{equation} J_{\nu}\left(z\right)\sim\frac{1}{\sqrt{2\pi\nu}}\left(\frac{ez}{2\nu}\right)^{\nu},\label{eq:asympJ}
\end{equation}
\begin{equation}
	Y_{\nu}\left(z\right)\sim-i{H^{(1)}_{\nu}}\left(z\right)\sim i{H^{(2)}_{\nu}}
\left(z\right)\sim-\sqrt{\frac{2}{\pi\nu}}\left(\frac{ez}{2\nu}\right)^{-\nu}.
\end{equation}
Therefore for the free-space Green's function, 
$$ |J_{p+1}(k_0 r) H_{p+1}(k_0 \rho)|/|J_{p}(k_0 r) H_{p}(k_0 \rho)| \sim \frac{r}{\rho}, $$
and the truncation errors of both the free-space multipole and local expansions decay exponentially 
as $p \to \infty$ \cite{rokhlin1990rapid}.

\vspace{0.1in}
{\bf \noindent Asymptotic Approximation of Integrals.} The asymptotic expansion for the Bessel functions can
be derived from the integral representations of these special functions. Note that these integral representations
are similar to the layered media Green's functions, therefore the same asymptotic analysis techniques can be 
applied, and the results can be used to derive more precise error bounds when truncating the layered media
multipole and local expansions. We demonstrate the ideas using the following examples.

We start with the integral representation of the Bessel function \cite{abramowitz1966handbook},
\begin{eqnarray*}
	J_n(z) & = & \frac{i^{-n}}{2 \pi} \int_{-\pi}^{\pi} e^{iz \cos \theta} \cos(n \theta) d \theta \\
	& = &  \frac{i^{-n}}{2 \pi} \int_{-\pi}^{\pi} \sum_{k=0}^{\infty} \frac{(iz \cos \theta)^k}{k!} \cos(n \theta) d \theta \\
	& = &  i\sum_{k=0}^{\infty} \frac{i^{-n}}{2 \pi} \frac{(iz)^k}{k!} \int_{-\pi}^{\pi}  (\cos \theta)^k \cos(n \theta) d \theta. 
\end{eqnarray*}
Note that when $k<n$, the integral $\int_{-\pi}^{\pi}  (\cos \theta)^k \cos(n \theta) d \theta =0$. By applying 
the Stirling's formula $ n! \sim \sqrt{2 \pi n} (\frac{n}{e})^n$ to the leading order of the expansion, one recovers the asymptotic 
form of the Bessel function for large orders in Eq.~(\ref{eq:asympJ}).
This asymptotic expansion is valid for all $z$ values and provides more accurate estimate when 
studying the truncation error in the multipole and local expansions, especially in the low frequency 
regime when $k_0 r$ and $k r$ are small.

For the layered media Green's function, both integrals in the basis functions $\Phi_p$ and $\Psi_p$ can be 
formulated as 
\begin{equation}
\int_{-\infty}^{\infty}  
    e^{-\sqrt{\lambda^{2}-k_1^2}y} e^{i\lambda x}  
    \left( \frac{\lambda - \sqrt{\lambda^2-k_2^2}}{k_2} \right)^p  
    \frac{\sigma (\lambda)}{4 \pi \sqrt{\lambda^{2}-k_1^2}}  d \lambda
\end{equation}
where $k_1$ is the wave number in the target layer and $k_2$ is the wave number in the source layer. The integral can be divided into three parts, 
when (a) $|\lambda| < \min\{k_1,k_2\}$, (b) $|\lambda|>\max\{k_1,k_2\}$, 
and (c) $\min\{k_1,k_2\} \leq |\lambda| \leq \max \{k_1,k_2\}$. The asymptotic expansion for each part can be 
derived using existing asymptotic analysis techniques for integrals \cite{wong2001asymptotic}. When $k_1=k_2$,
the integral in (a) is often referred to as the ``propagating" part and the integral in (c) becomes the ``evanescent" 
part. The asymptotic properties of the basis functions $\Phi_p$ and $\Psi_p$ are determined by the ``propagating" part 
for large $|| \mathbf{x}||$ values with fixed $p$, and by the ``evanescent" part for large $p$ values with
fixed $||\mathbf{x}||$. To understand the truncation errors in the multipole and local expansions of the layered media
Green's function, we therefore focus on the evanescent part and demonstrate the asymptotic analysis for the simplified 
integral
\begin{equation}
	\int_{0}^{\infty}  
	e^{- t y} e^{i \sqrt{t^2+k_1^2}  x}  
	\left( \frac{\sqrt{t^2+k_1^2} + \sqrt{t^2+k_1^2-k_2^2}}{k_2} \right)^p  
	\frac{\tilde{\sigma}(t)}{4 \pi \sqrt{t^{2}+k_1^2}}  dt 
	\label{eq:asympLayer}
\end{equation}	
where we assume $k_1>k_2$. The case when $k_1<k_2$ can be analyzed in a similar way.
Instead of deriving the asymptotic expansion of Eq.~(\ref{eq:asympLayer}) directly, 
we adopt the following steps to further simplify the integral to a more standard form that is commonly
used when analyzing the asymptotic behavior of the Bessel functions. First, using the polar coordinates
$(r, \theta)$ of $(x,y)$, we rewrite the $\{x,y\}$-related exponential part as
$$ e^{- t y} e^{i \sqrt{t^2+k_1^2} x}  =  e^{-k_1 r \left( \frac{t \sin (\theta )}{k_1}-\frac{i \cos (\theta ) \sqrt{k_1^2+t^2}}{k_1} \right)}. $$  
Second, we define a new variable 
$u=\frac{t \sin (\theta )}{k_1}-\frac{i \cos (\theta ) \sqrt{k_1^2+t^2}}{k_1} $. 
Clearly, the new integral for the $u$-variable is on a complex contour, not on the real axis.
Third, using contour integration and residue theorem, we can rewrite the integral on the complex contour
back to the sum of an integral on the real axis and an easy-to-analyze integral on a line segment in the complex plane 
(which will be explained in detail in Section~\ref{sec:altdirection}). Finally, we focus on the asymptotic 
expansion of the dominating integral
\begin{eqnarray*}
\int_{0}^{\infty}  
    {e^{- k_1 r u} } 
    {\left(  \frac{ \beta(u) + \sqrt{\beta^2(u)-k_2^2}}{k_2} \right)^p } 
    {\frac{\dbtilde{\sigma} (u)}{4 \pi \sqrt{1+u^2} }  } d u
= \int_{0}^{\infty}  e^{- k_1 r u} \left(f(u)\right)^p g(u) du
\end{eqnarray*}	
where $\beta(u)=k_1 \left(\sqrt{u^2+1} \sin (\theta )-i u \cos (\theta )\right)$ and $\dbtilde{\sigma} (u) = \tilde{\sigma}(t(u))$.
For any fixed $z=k_1 r$, one approach to find the asymptotic behavior of this integral is to 
first perform another change of variable $\lambda = k_1 r u=z u$ (or $u=\frac{\lambda}{k_1 r} = \frac{\lambda}{z}$)
and consider the new integral 
$$ \int_{0}^{\infty}  e^{- \lambda } \left (f( \frac{\lambda}{z})\right)^p g( \frac{\lambda}{z}) \frac{d \lambda}{z}
 = (k_2 r)^{-p} \int_{0}^{\infty}  e^{- \lambda }\left(   \tilde{\beta}(z) + \sqrt{ \tilde{\beta}^2(z)-(\frac{k_2}{k_1} z)^2} \right)^p  
	\tilde{g}( \frac{\lambda}{z}) d \lambda  $$
where $\tilde{\beta}(z)=\sqrt{\lambda ^2+z^2} \sin (\theta )  -i \lambda  \cos (\theta ).$ 
Under proper conditions, the function $ \left(   \tilde{\beta}(z) + \sqrt{ \tilde{\beta}^2(z)-(\frac{k_2}{k_1} z)^2} \right)^p  
\tilde{g}( \frac{\lambda}{z})$ is an analytic function for $z$ values on the right half complex plane away from the
origin, we can consider its Taylor expansion 
$$\left(   \tilde{\beta}(z) + \sqrt{ \tilde{\beta}^2(z)-(\frac{k_2}{k_1} z)^2} \right)^p  
\tilde{g}( \frac{\lambda}{z}) = \sum_{k=0}^{\infty} t_k (\lambda,p,\frac{k_2}{k_1}) z^k $$
as a function of $z$. Then, the integral representation of evanescent part can be derived as
$$  (k_2 r)^{-p} \sum_{k=0}^{\infty}  \left( \int_{0}^{\infty}  e^{- \lambda } t_k (\lambda,p,\frac{k_2}{k_1}) d\lambda \right) z^k.
$$
The leading order term of the expansion for very large $p$ values is approximately 
$$ (k_2 r)^{-p} {e^{i p \theta_0}} \int_{0}^{\infty}  e^{- \lambda } (2 \lambda)^p  d\lambda$$
for some constant $\theta_0$, which has the same asymptotic properties as the Hankel function $H_p(k_2 r) e^{i p \theta_0}.$
Note that the expansion can be used to study the properties of both the cases when $p \to \pm \infty$ or when $z$ 
is small (low-frequency regime). Without presenting the details, we summarize our results in the following 
theorem.
\begin{theorem}
{\bf (1)} The multipole and local expansions satisfy the following convergence estimates for 
	large $p$ values,  
\begin{align}
& |J_{p+1}(k_1 r) \Phi_{p+1}(x,y)|/|J_{p}(k_1 r) \Phi_{p}(x,y)| \sim \frac{r}{\rho}, \\
& |J_{p+1}(k_2 r) \Psi_{p+1}(x_0,y_0)|/|J_{p}(k_2 r) \Psi_{p}(x_0,y_0)| \sim \frac{r}{\tilde{\rho}}
\end{align}
	where $r$ is the distance from the particle to its box center and $\rho$ (or $\tilde{\rho}$) is the ``modified" 
	distance between the box center and far-field point. Both $d$ and $\pm$ sign are
	considered to correct the $\rho$ (or $\tilde{\rho}$) value. For example, $\rho=\sqrt{(y-y_c^s+d)^2+(x-x_c^s)^2}$ 
	for $\Phi(x,y)$ in Eq.~(\ref{eq:mpbasis}) and $\tilde{\rho}=\sqrt{ (y_c^t+d-y_0)^2+(x_c^t-x_0)^2}$
	for $\Psi(x_0,y_0)$ in Eq.~(\ref{eq:localbasis}). \\
{\bf (2)} When the wave numbers $k_1$ and $k_2$ are small, the required number of terms in the 
multipole and local expansions for a prescribed accuracy in the layered media FMM is approximately the same as that 
in the free-space FMM.
\end{theorem}
This theorem simply states that both the multipole and local expansions in the layered media FMM are 
exponentially convergent and presents the asymptotic convergence rates. In the low-frequency regime, 
it provides more precise estimates of the number of terms required.

{\bf \noindent Remark:} Note that 
using the ``modified" distance, the truncated multipole expansion of the scattered field Green's 
function may become valid for a neighboring (including self) target box, and therefore can be translated to its neighbor's 
local expansion using the M2L operator instead of the more expensive S2T operator. 

\vspace{0.1in}
{\bf \noindent Laplace Transform and Complex Images.} For many layered media Green's functions,
one can justify that the inverse Laplace transform of $\dbtilde{\sigma}(z)$
can be derived using the Fourier-Mellin integral formula
$$ \mathcal{L}^{-1} \{\dbtilde{\sigma} \}(t) = \frac{1}{2 \pi i} 
 \lim_{T \to \infty} \int_{\gamma -iT}^{\gamma + iT} e^{zt} \dbtilde{\sigma} (z) d z  $$ 
where the integration is done along the vertical line $\text{Re}(z)=\gamma$ in the complex plane
and $\gamma$ is greater than the real part of all singularities of $ \dbtilde{\sigma} (z)$,
so that the image part in the layered media Green's function can be represented as
$$ \dbtilde{\sigma} (\lambda) = \int_0^{\infty} e^{-\sqrt{\lambda^2-k_0^2} (t - \gamma_0)  } f(t) dt$$
for some $\gamma_0$ which may or may not be the same as $\gamma$.
Plugging this representation in the original layered media Green's function in Eq.~(\ref{eq:Green}), we
have 
\begin{equation}
G(\mathbf{x}, \mathbf{x}_0)= \int_0^{\infty} \hskip-0.1in \int_{-\infty}^{\infty} 
    \frac{ e^{-\sqrt{\lambda^{2}-k^2}(y+d)} e^{i\lambda x} }{ {4 \pi \sqrt{\lambda^{2}-k^2}}}
    { e^{ \pm \sqrt{\lambda^2-k_0^2} y_0 - \sqrt{\lambda^2-k_0^2} (t - \gamma_0) }e^{-i\lambda x_0 } } 
    {f(t)} d\lambda dt.
\end{equation}
In this new representation, $f(t)$ can be considered as the complex image located at $(x_0,\pm y_0 -t)$.
For many layered media, it is sufficient to analyze the convergence of the multipole and local
expansions of the layered media Green's functions for {\it each} $t$-mode where the new source box center is 
at $(x_c^s, \pm y_c^s -(t-\gamma_0))$. In \cite{cho2018heterogeneous}, 
this complex image approach is applied to the 2-D half-space layered medium Green's function with impedance 
boundary conditions (see Example 3 of Sec. \ref{sec:Green}), and it becomes straightforward to verify using the complex
images that for the same prescribed accuracy requirement, the number of expansion terms for the layered medium 
case is no more than that for the free-space case.  Another application of the complex image approach is 
for the direct source-target (S2T) interactions. For a neighboring source box, when all the complex images 
of the scattered field are well-separated from the target box, the source box's multipole expansion becomes 
applicable, and it is more efficient to translate the multipole expansion to a local expansion of the 
target box. This technique is used in  \cite{cho2018heterogeneous} to compute the scattered field 
part of the source-target interactions.

\vspace{0.1in}
{\bf \noindent Precomputed Tables for Number of Expansion Terms.} 
In practice, both the asymptotic expansion and complex image approach only give a very rough estimate of the number of terms
required in different expansions of the layered media Green's function. A more practical approach is to precompute
a table (or table of tables) for different layered media settings. This approach is problem dependent. We are constructing such tables for several real world applications, and results will be reported in the future.

\subsection{Alternative Direction Sommerfeld Integral Representations}
\label{sec:altdirection}
Another difficulty in simulating waves in layered media is when the source and target are very
close to each other or when the source and target are close to the interface of different layers
for the scattered field. In these cases, the computation of layered media Green's function and M2L translation operator becomes extremely expensive.  For example, when $y-y_0>0$ is close to zero and $|x-x_0|$ is relatively a big number in the Sommerfeld integral representation of free-space Green's function in Eq.~(\ref{eq:spectralfree}) or the half-space layered medium Green's function with impedance boundary conditions in Eq.~(\ref{eq:halfspace}),  the exponential term $e^{-\sqrt{\lambda^{2}-k^2}(y-y_0)}$ decays slowly and the highly oscillatory term $e^{i\lambda (x-x_0) }$ must be sufficiently sampled. A similar problem occurs in the M2L operator in Eq.~(\ref{eq:M2L}). When the ratio of $(y_c^t-y_c^s+d)/|x_c^t-x_c^s|$ is small, a wide range of $\lambda$ values have to be sampled before the integrand decays to zero sufficiently.

The numerical difficulty is not from any inherent properties of the original physical problem, 
it is the result of using the inefficient integral representation in the numerical computation. 
For the same Hankel function $H_0(\beta r)$, which is the Green's function of the free-space
Helmholtz equation, a common practice is to divide $\ensuremath{{\bf R}}^2$ into four
{\it overlapping regions}--North, South, East, West--corresponding to points
$(x,y)\in \ensuremath{{\bf R}}^2$ with $y>0,y<0,x>0,x<0$, respectively. In each region, the plane
wave representation of $H_0(\beta r)$ takes the following forms:
\begin{align}
&H_0(\beta r) = \\
&\left\{\begin{array}{l}
\frac{1}{\pi} \int_{0}^{\pi}  e^{i \beta(y\sin\theta-x\cos\theta)} d \theta +  
 \frac{1}{i\pi}\int_0^{\infty} \frac{e^{-ty}}{\rho_{\beta}(t)}\left( e^{i\rho_{\beta}(t)x}+e^{-i\rho_{\beta}(t)x}\right) dt  
 ~~\mbox{  \textbf{North},} \\
\frac{1}{\pi} \int_{0}^{\pi}  e^{i\beta(-y\sin\theta-x\cos\theta)} d \theta +  
 \frac{1}{i\pi}\int_0^{\infty} \frac{e^{ ty}}{\rho_{\beta}(t)}\left( e^{i\rho_{\beta}(t)x}+e^{-i\rho_{\beta}(t)x}\right) dt  
 \mbox{ \textbf{South},} \\
\frac{1}{\pi} \int_{0}^{\pi}  e^{i\beta(-y\cos\theta +x\sin\theta)} d \theta +  
 \frac{1}{i\pi}\int_0^{\infty} \frac{e^{-tx}}{\rho_{\beta}(t)}\left( e^{i\rho_{\beta}(t)y}+e^{-i\rho_{\beta}(t)y}\right) dt  
 \mbox{ \textbf{East},}  \\
\frac{1}{\pi} \int_{0}^{\pi}  e^{i\beta(-y\cos\theta -x\sin\theta)} d \theta +  
 \frac{1}{i\pi}\int_0^{\infty} \frac{e^{ tx}}{\rho_{\beta}(t)}\left( e^{i\rho_{\beta}(t)y}+e^{-i\rho_{\beta}(t)y}\right) dt  
 \mbox{ \textbf{West}}
\end{array}\right. \nonumber
\end{align}
where $\rho_{\beta}(t)=\sqrt{t^2+\beta^2}$ and  $(r, \theta)$ are the polar coordinates of the point $(x,y)$. 
For higher order Hankel functions $H_l(\beta r)e^{i l \theta}$, we have the following integral representations.  
\begin{align}
&H_l(\beta r)e^{i l \theta} = \nonumber\\
&\left\{\begin{array}{l}
\frac{i^l}{\pi} \int_{0}^{\pi} e^{i\beta( y\sin\theta - x\cos\theta) }e^{-i l \theta} d\theta + \nonumber\\ 
\frac{(-i)^l}{i \pi}\int_0^{\infty} \frac{e^{-ty}}{\rho_{\beta}(t)} 
     \left( e^{i \rho_{\beta}(t) x} \left( \frac{ \rho_{\beta}(t)-t}{\beta}\right)^l+ \nonumber
         e^{-i\rho_{\beta}(t) x} \left( \frac{-\rho_{\beta}(t)-t}{\beta}\right)^l \right) dt \mbox{ \textbf{North},}\\
\frac{i^l}{\pi} \int_{0}^{\pi} e^{i\beta(-y\sin\theta - x\cos\theta) }e^{ i l \theta} d\theta + \nonumber\\ 
\frac{(-i)^l}{i \pi}\int_0^{\infty} \frac{e^{ty}}{\rho_{\beta}(t)} 
     \left( e^{i \rho_{\beta}(t) x} \left( \frac{ \rho_{\beta}(t)+t}{\beta}\right)^l+ 
         e^{-i\rho_{\beta}(t) x} \left( \frac{-\rho_{\beta}(t)+t}{\beta}\right)^l \right) dt \mbox{ \textbf{South},}\\
\frac{ 1}{\pi} \int_{0}^{\pi} e^{i\beta( x\sin\theta - y\cos\theta) }e^{ i l \theta} d\theta + \\ 
\frac{(-1)^l}{i \pi}\int_0^{\infty} \frac{e^{-tx}}{\rho_{\beta}(t)} 
     \left( e^{i \rho_{\beta}(t) y} \left( \frac{-\rho_{\beta}(t)-t}{\beta}\right)^l+ 
         e^{-i\rho_{\beta}(t) y} \left( \frac{ \rho_{\beta}(t)-t}{\beta}\right)^l \right) dt \mbox{ \textbf{East},}\\
\frac{(-1)^l}{\pi} \int_{0}^{\pi} e^{i\beta(-x\sin\theta - y\cos\theta) }e^{-i l \theta} d\theta + \\ 
\frac{(-1)^l}{i \pi}\int_0^{\infty} \frac{e^{ tx}}{\rho_{\beta}(t)} 
     \left( e^{i \rho_{\beta}(t) y} \left( \frac{-\rho_{\beta}(t)+t}{\beta}\right)^l+ 
         e^{-i\rho_{\beta}(t) y} \left( \frac{ \rho_{\beta}(t)+t}{\beta}\right)^l \right) dt \mbox{ \textbf{West}}
 \end{array}\right.
\end{align} 
where the first and the second integrals in each formula are called the propagating part and evanescent part, respectively.
These directional representations have been applied in the {\it new version} of FMM
in both two and three dimensions \cite{crutchfield2006remarks,greengard1997new,cheng2006wideband,greengard2002new}.
They are also effective tools to compute the lattice sums of the free-space Green's functions 
\cite{dienstfrey2001lattice,huang1999integral}. Clearly in these formulas, the oscillation of the integrand 
in the propagating part is controlled by $\beta r$ and hence no numerical quadrature issues arise. 
For the evanescent part at the overlapping regions, e.g., when both $x,y>0$, both the 
north and east integral representations are valid and can be applied, however their numerical properties 
are very different. Clearly when $y/x\gg1$, the north formula is preferred, and when $y/x\ll 1$, the east formula
can be computed more efficiently using existing quadrature techniques. 
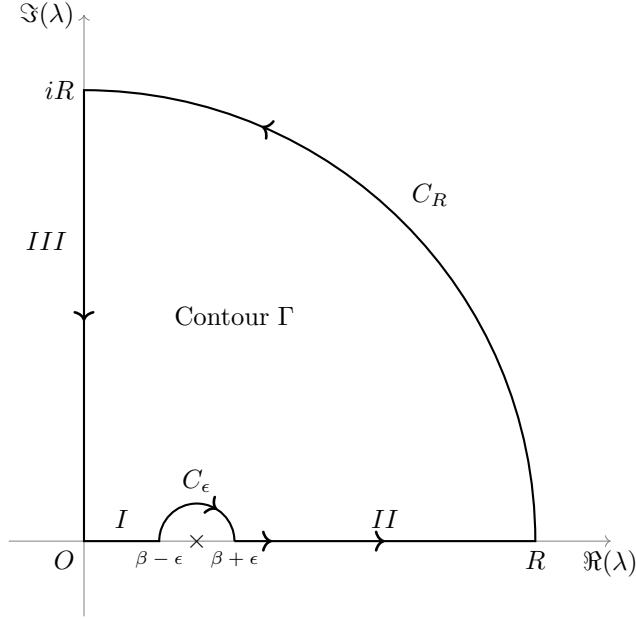
\begin{figure}[t]
  \centering
  \begin{tikzpicture} 
[
    decoration={%
      markings,
      mark=at position 0.5cm with {\arrow[line width=1pt]{>}},
      mark=at position 2cm with {\arrow[line width=1pt]{>}},
      mark=at position 0.5 with {\arrow[line width=1pt]{>}},
      mark=at position 0.75 with {\arrow[line width=1pt]{>}},
      mark=at position -5mm with {\arrow[line width=1pt]{>}},
    }
  ]
  \draw [help lines,->] (-1,0) -- (7,0) coordinate (xaxis);
  \draw [help lines,->] (0,-1) -- (0,7) coordinate (yaxis);
  \node at (1.5,0) {$\times$};
  \path [draw, line width=0.8pt, postaction=decorate] (2,0) node [below, font=\scriptsize] {$\beta+\epsilon$} -- 
    (6,0) node [below] {$R$} arc (0:90:6) node [left] {$iR$} -- (0,0) -- 
    (1,0) node [below, font=\scriptsize] {$\beta-\epsilon$} arc (180:0:.5);
  \node [below] at (xaxis) {$\Re(\lambda)$};
  \node [left] at (yaxis) {$\Im(\lambda)$};
  \node [below left] {$O$};
  \node at (1.5,.8) {$C_{\epsilon}$};
  \node at (2,3) {Contour $\Gamma$};
  \node at (4.6,4.6) {$C_{R}$};
  \node at (0.5,0.3) {$I$};
  \node at (4,0.3) {$II$};
  \node at (-0.5,4) {$III$};

\end{tikzpicture} 
  \caption{Contour Integral}
  \label{fig:contour}
\end{figure}

The alternative direction integral can be derived using contour integration and residue theorem in 
two dimensions. To demonstrate how this technique works, we consider the contour integral for $x,y>0$ 
in fig.~\ref{fig:contour}
\begin{equation}
 \frac{1}{\pi} \int_{\Gamma} \frac{e^{i(\lambda x + \sqrt{\beta^2-\lambda^2} y)}}{\sqrt{\beta^2-\lambda^2}} d \lambda. 
\end{equation}
Because there are no singularities or branch cut points inside the contour $\Gamma$, we have 
\begin{equation} 
 \frac{1}{\pi} \int_{\Gamma} \frac{e^{i(\lambda x + \sqrt{\beta^2-\lambda^2} y)}}{\sqrt{\beta^2-\lambda^2}} d \lambda = \int_I + \int_{C_{\epsilon}} + \int_{II} + \int_{C_R} + \int_{III}=0. 
\label{eq:alternative} 
\end{equation}
When $\epsilon \to 0$ and $R \to \infty$,
	\begin{align}
	&\frac{1}{\pi} \int_{C_{\epsilon}} \frac{e^{i(\lambda x + \sqrt{\beta^2-\lambda^2} y)}}{\sqrt{\beta^2-\lambda^2}} d \lambda  \to  0~~ ,~~ \frac{1}{\pi}\int_{C_R}\frac{e^{i(\lambda x + \sqrt{\beta^2-\lambda^2} y)}}{\sqrt{\beta^2-\lambda^2}} d \lambda  \to 0 \nonumber
\end{align}
and
	\begin{align}
	&\frac{1}{\pi}\int_I  \frac{e^{i(\lambda x + \sqrt{\beta^2-\lambda^2} y)}}{\sqrt{\beta^2-\lambda^2}} d \lambda \to  \frac{1}{\pi} \int_0^{\beta} \frac{e^{i (\lambda x -\sqrt{\beta^2-\lambda^2}y)}}{-\sqrt{\beta^2-\lambda^2}} d \lambda,\nonumber\\
	&\frac{1}{\pi}\int_{II}\frac{e^{i(\lambda x + \sqrt{\beta^2-\lambda^2} y)}}{\sqrt{\beta^2-\lambda^2}} d \lambda  \to   \frac{1}{\pi} \int_{\beta}^{\infty} \frac{e^{i (\lambda x + i \sqrt{\lambda^2-\beta^2}y)}}{i \sqrt{\lambda^2 - \beta^2}} d \lambda, \nonumber\\
	&\frac{1}{\pi}\int_{III} \frac{e^{i(\lambda x + \sqrt{\beta^2-\lambda^2} y)}}{\sqrt{\beta^2-\lambda^2}} d \lambda  \to  \frac{1}{\pi} \int_{0}^{\infty} \frac{i e^{ (-t x - i \sqrt{t^2+\beta^2}y})}{ \sqrt{t^2 + \beta^2}} dt.  \nonumber 
	\end{align}
Therefore, the evanescent part in $\int_{II}$ can be computed using $-\left( \int_{I} + \int_{III} \right)$. In fact, the east plane wave representation can be re-derived from the north one using this approach for $x,y>0$.

Finally, we point out that the proper choice of the integration contour is problem dependent. We present more details in 
Section~\ref{sec:numerical} for other layered media Green's functions. 

\section{Preliminary Numerical Experiments} 
\label{sec:numerical}
In this section, we present numerical experiments to validate our theoretical analysis of the 
general numerical framework. Matlab and Mathematica codes are developed to numerically validate the analysis presented in Section \ref{sec:analysis}. 

\vspace{0.1in}
{\noindent \bf Alternative Direction Sommerfeld Integral Representation.}
We have studied and validated the alternative direction Sommerfeld integral formulas for several layered media Green's 
functions using Mathematica. Different direction plane wave representations of the free-space Green's function presented 
in Section~\ref{sec:altdirection} can be readily found from existing literature 
\cite{dienstfrey2001lattice,huang1999integral}, therefore we focus on the results for the half-space 
layered medium with impedance boundary condition. The three layered medium Green's functions can be handled in a very similar way, we
skip the detailed formulas and interested readers are referred to the Mathematica files for these formulas as well as their validations.

For the half-space layered medium Green's function, we focus on the evanescent part of the Green's function given by 
\begin{equation} 
\label{eq:evan1}
\int_0^{\infty} \frac{e^{-t y} e^{i x \sqrt{t^2+k^2}}}{\sqrt{t^2+k^2}} \frac{(t+i \alpha)}{(t-i \alpha)} dt.
\end{equation}
To avoid the pole on the imaginary axis, we have numerically tested the following two contours.
\begin{figure}[t]
  \centering
  \begin{tikzpicture} 
[
    decoration={%
      markings,
      mark=at position 0.5cm with {\arrow[line width=1pt]{>}},
      mark=at position 0.5 with {\arrow[line width=1pt]{>}},
      mark=at position 0.75 with {\arrow[line width=1pt]{>}},
    }
  ]
  \draw [help lines,->] (-1,0) -- (4.2,0) coordinate (xaxis);
  \draw [help lines,->] (0,-1) -- (0,4) coordinate (yaxis);
  \node at (1,0) {$\times$};
  \path [draw, line width=0.8pt, postaction=decorate] (1,0) node [below, font=\scriptsize] {$c$} -- 
    (3.5,0) node [below] {$R$} arc (0:73.4:3.5) -- (1,0);  
  \node [below] at (xaxis) {$\Re(\lambda)$};
  \node [left] at (yaxis) {$\Im(\lambda)$};
  \node [below left] {$O$};
  \node at (2.0,3.9) {$1^{st}$ contour};
  \node at (1.9,0.3) {$I$};
  \node at (0.5,-0.3) {$IV$};
  \node at (3,2.3) {$II$};
  \node at (0.5,1.9) {$III$};

  \draw [help lines,->] (5.2,3.5) -- (11,3.5) coordinate (xaxis);
  \draw [help lines,->] (6,-1) -- (6,4) coordinate (yaxis);
  \draw[->,line width=0.8pt,domain=0:2,smooth,variable=\t,red]  plot ({(\t)*0.7+6},{3-(sqrt(\t*\t+1)*0.714});
  \draw[line width=0.8pt,domain=2:4,smooth,variable=\t,red]  plot ({(\t)*0.7+6},{3-(sqrt(\t*\t+1)*0.714});
  \draw[->,line width=0.8pt,domain=-50.89:-25,smooth,variable=\t,red]  plot ({4.44*cos(\t)+6},{3.5+4.44*sin(\t)});
  \draw[line width=0.8pt,domain=-25:0,smooth,variable=\t,red]  plot ({4.44*cos(\t)+6},{3.5+4.44*sin(\t)});
  \draw[->,line width=0.8pt,red] (10.44,3.5) -- (8.7,3.5) ;
  \draw[line width=0.8pt,red] (8.7,3.5) -- (7.6,3.5) ;
  \draw[->,line width=0.8pt,red] (7.6,3.5) -- (6.8, 2.9);
  \draw[line width=0.8pt,red](6.8, 2.9) -- (6,2.29);
  \node [below] at (xaxis) {$\Re(\lambda)$};
  \node [left] at (yaxis) {$\Im(\lambda)$};
  \node [below left] {$O$};
  \node at (9.0,3.9) {$2^{nd}$ contour};
  \node at (7.3,0.9) {$I$};
  \node at (10.5,1.3) {$II$};
  \node at (8.7,3.3) {$III$};
  \node at (7.2,2.9) {$IV$};
  \node at (10.45,3.8) {$R$};
  \node at (7.6,3.7) {$c$};

\end{tikzpicture} 
  \caption{Two different contours for half-space layered medium with impedance boundary condition.}
  \label{fig:contour2}
\end{figure}
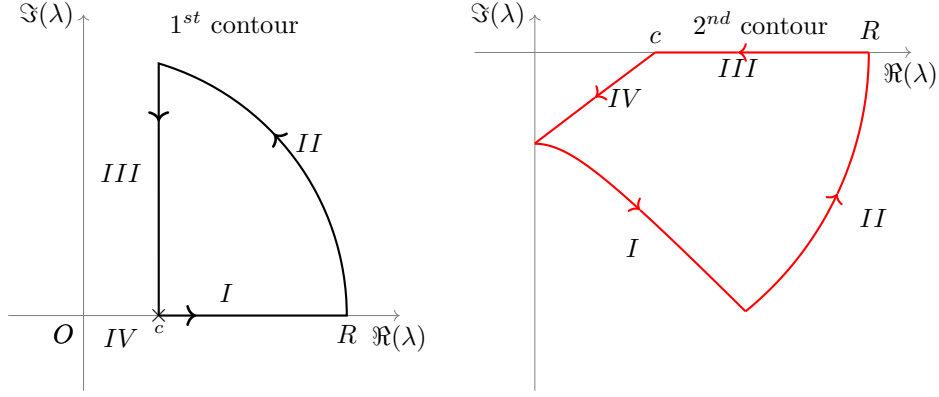
In the first contour (left of fig.~\ref{fig:contour2}), a positive $c$ value is chosen so that
\[ 
\int_0^{\infty} \frac{e^{-t y} e^{i x \sqrt{t^2+k^2}}}{\sqrt{t^2+k^2}} \frac{(t+i \alpha)}{(t-i \alpha)} dt
= \left( \int_0^{c} + \int_c^{\infty} \right) \frac{e^{-t y} e^{i x \sqrt{t^2+k^2}}}{\sqrt{t^2+k^2}} \frac{(t+i \alpha)}{(t-i \alpha)} dt. 
\]
As the first term integrates from $0$ to $c$ on a finite (and reasonably small by proper choice of $c$) line segment 
(labeled $IV$), it can therefore be efficiently evaluated using standard Gauss quadrature rules. For the second integral, 
as the sum of the contour integrals on $I+II+III$ is $0$ and the integral on $II$ approaches $0$ when $R \to \infty$, 
the alternative direction representation of the integral $$ \int_c^{\infty} \frac{e^{-t y} 
e^{i x \sqrt{t^2+k^2}}}{\sqrt{t^2+k^2}} \frac{(t+i \alpha)}{(t-i \alpha)} dt $$  is given by
\begin{equation}
\label{eq:evan2}
i e^{-cy} \int_0^{\infty} \frac{e^{- \lambda x } e^{-i \lambda y} }{\sqrt{(c+i\lambda)^2+k^2}} 
e^{ix \left( \sqrt{(c+i \lambda)^2+k^2} - \sqrt{(i \lambda)^2}\right)} \frac{((c+i\lambda)+i \alpha)}{((c+i \lambda)-i \alpha)} d\lambda. 
\end{equation}
This representation is numerically validated using Mathematica's {\tt{NIntegrate}} with options 
{\tt{AccuracyGoal}$\to$20}, {\tt{PrecisionGoal}$\to$20}, {\tt{WorkingPrecision}$\to$60}, {\tt{MaxRecursion}$\to$100},
{\tt{Method}$\to$\tt{DoubleExponential}} for different $x,y>0$ values.

In the second contour (right of fig.~\ref{fig:contour2}), we assume $(r,\theta)$ are the polar coordinates 
of $(x,y)$, and perform the change of variable
$u=\frac{t}{k} \sin \theta - i \frac{\sqrt{t^2+k^2}}{k} \cos \theta$, then the evanescent part becomes
\[ \int_{I} \frac{e^{-k r u}}{\sqrt{1+u^2}} \frac{ i k \cos \theta \sqrt{1+u^2} + u k \sin \theta + i \alpha}{ i k 
  \cos \theta \sqrt{1+u^2} + u k \sin \theta - i \alpha} du
\]
where the contour $I$ is a curve defined by $z(t)=  \frac{t}{k} \sin \theta - i \frac{\sqrt{t^2+k^2}}{k}$, $t \in [0, \infty)$.
As the integral on $II$ approaches $0$ when $R \to \infty$, the evanescent part becomes the (negative) sum of the 
integrals on $III$ and $IV$ given by
\begin{equation}
 \int_{c}^{\infty} \frac{e^{-k r u}}{\sqrt{1+u^2}} \frac{ \phi(u)+ i \alpha}{\phi(u)- i \alpha} du 
\label{eq:evan3}
\end{equation}
and 
\begin{equation}
	\label{eq:evan3part2}
 \int_{0}^{1} \frac{e^{-k r \psi(u)}}{\sqrt{1+\psi(u)^2}} \frac{ \phi(\psi(u))+ i \alpha}{\phi(\psi(u))- i \alpha} (c+ i \cos \theta) du
\end{equation}
where $ \phi(u)= i k \cos \theta \sqrt{1+u^2} + u k \sin \theta$, $\psi(u)= (c+i \cos \theta) u-i \cos \theta$, and $c$
is a constant to be optimized so that both the Gauss quadrature applied on $IV$ and Laguerre quadrature on $III$ 
converge rapidly. These formulas are also validated using Mathematica for different $(x,y)$ and $c$ values. 
Note that in the second contour, the oscillatory term $e^{i \lambda x}$ or $e^{i \lambda y}$ is completely removed from the 
integrand, at the cost of a new rational function on a contour closer to the singularities. 

\begin{figure}[t]
   \centering
	\includegraphics[width=4.8in]{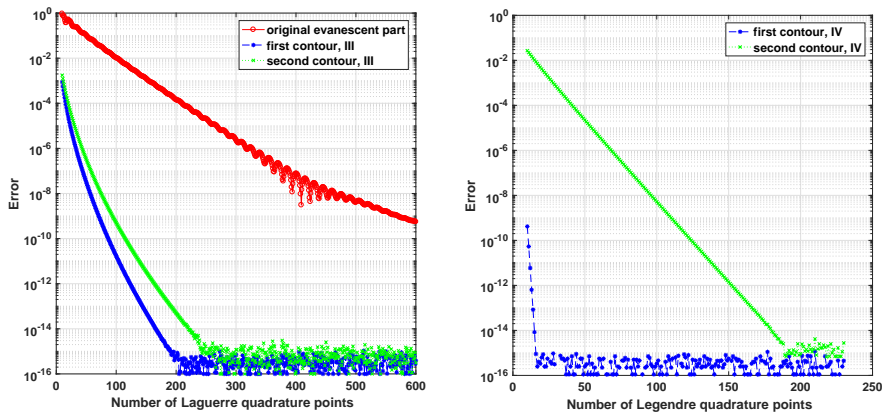}
   \caption{Convergence of the quadrature rules for different integrals.}
  \label{fig:quadconv}
\end{figure}
We have tested numerically the convergence of the quadrature rules for the original and mathematically equivalent 
alternative direction representations of the half-space layered medium Green's function. In fig.~\ref{fig:quadconv}, 
we set $x=1$, $y=0.1$, $k=1$, and $c=2$, and present the numerical errors when different numbers of nodes are used
in the quadrature rules for different integrals. The reference solutions are computed using Mathematica requesting 
more than $20$ correct digits. On the left of fig.~\ref{fig:quadconv},  we present the accuracy when different 
numbers of Laguerre quadrature nodes are used for the original evanescent part in Eq.~(\ref{eq:evan1}), the integral 
on $III$ of the first contour in Eq.~(\ref{eq:evan2}), and the integral on $III$ of the second contours in 
Eq.~(\ref{eq:evan3}). As $y\ll x$, the original integral converges slowly due to the oscillatory term $e^{i \sqrt{t^2+k^2} x}$. 
The Laguerre quadrature for the alternative direction integrals, on the other hand, converges much faster. 
For a fair comparison, we also present
the convergence of the quadrature rules for the integrals on the finite line segments. When Legendre polynomial based
Gauss quadrature is applied to the integral $IV$ on both the first and second contours, for double precision requirement,
under the current setting of $c$, approximately $17$ nodes are required for the integral $IV$ on the first integral, 
while around $200$ Gauss nodes are required for the integral $IV$ on the second contour. In our numerical implementation, 
as the formulas for the first contour is easy to derive and manipulate, we therefore adopt the first contour for the 
evanescent part of a general layered media Green's function.
\begin{figure}[t]
   \centering
	\includegraphics[width=4.7in]{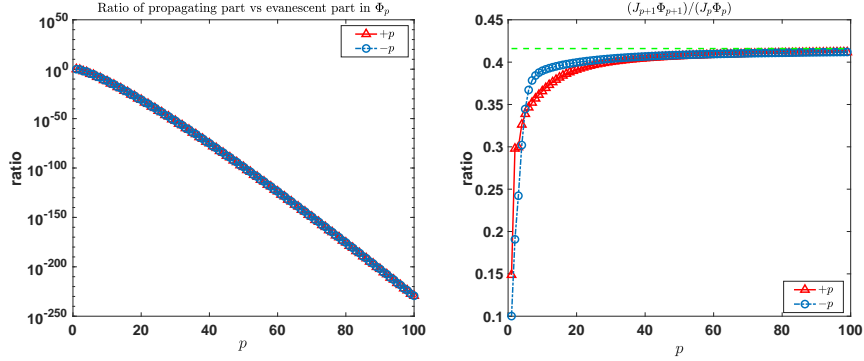}
   \caption{Convergence of the multipole expansions for impedance half-space layered medium.}
  \label{fig:orderconv}
\end{figure}
\vspace{0.1in}

{\noindent \bf Convergence of Multipole and Local Expansions.}
As we have discussed in Sec.~\ref{sec:truncate}, the convergence of the multipole and local expansions is determined by
the ratio of $r / \rho$, where $r$ is the distance between the source (target) and center of the source (target)
box in the multipole (local) expansion, and $\rho$ is the ``modified" distance between the far-field target (source) 
and center of the source (target) box (constant $d$ and $\pm$ sign are considered). We have numerically validated the 
analysis. In fig.~\ref{fig:orderconv}, we present the results for the half-space layered medium with impedance 
boundary conditions using settings $x_{target}-x_c^s=2$, $y_{target}+y_c^s=3$, $r=1.5$, $k=0.1$, $\alpha=1$, and
the modified distance $\rho=\sqrt{(x_{target}-x_c^s)^2+(y_{target}+y_c^s)^2}$.

\begin{figure}[t]
   \centering
	\includegraphics[width=3.8in]{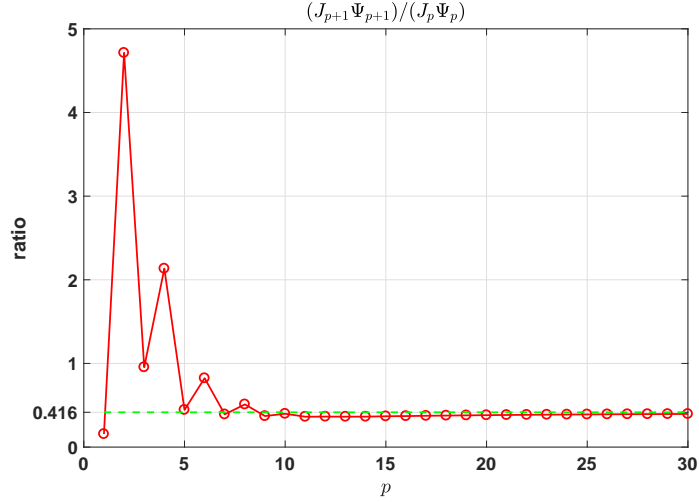}
   \caption{Convergence of the local expansion for three-layered medium Green's function.}
  \label{fig:order3layer}
\end{figure}
On the left of the figure, we compare the propagating and evanescent parts of the basis $\Phi_p$ for different $p$
values. Clearly, for large order $p$, the propagating part can be neglected when analyzing the truncation errors.
On the right of the figure, we plot the ratio $(J_{p+1} \Phi_{p+1}) / (J_p \Phi_p) $ when $p \to \pm \infty$. When
$|p|$ increases, the ratio approaches the constant $r / \rho$ (dashed green line).

A similar analysis is also performed for the three-layered medium Green's function with wave numbers $k_1$, $k_2$, and $k_3$ 
in each layer. In the numerical experiment, we set $k_1=1$, $k_2=3$, $k_3=1$ and consider the contribution from upper
layer ($y=0$) 
\[ g_2^t(\mathbf{x},\mathbf{x}_{0}) =  \int_{-\infty}^{\infty}
	{ e^{\sqrt{\lambda^2-k_2^{2}} y }e^{i\lambda x} }
	{e^{-\sqrt{\lambda^2-k_1^{2}} y_0 }e^{-i\lambda x_0} }
	{ \frac{1}{4\pi \sqrt{\lambda^2-k_2^2} } \sigma_2^t(\lambda)} d\lambda
\]
where \\
\scalebox{0.96}{ $\sigma_2^t(\lambda)=\frac{e^{d \sqrt{\lambda ^2-k_2^2}} (\lambda ^2+\sqrt{\lambda ^2-k_2^2} \sqrt{\lambda ^2-k_3^2}-k_2^2)}{\sinh (d \sqrt{\lambda ^2-k_2^2}) (\lambda ^2+\sqrt{\lambda ^2-k_1^2} \sqrt{\lambda ^2-k_3^2}-k_2^2)+\sqrt{\lambda ^2-k_2^2} (\sqrt{\lambda ^2-k_1^2}+\sqrt{\lambda ^2-k_3^2}) \cosh (d \sqrt{\lambda ^2-k_2^2})}. 
$ } \\
We study the convergence of the local expansion where the basis function $\Psi_p(x_0,y_0)$ is given by \\
\begin{align}
\Psi_p(x_0,y_0)= 
\int_{-\infty}^{\infty}  &
         {\left( \frac{\lambda - i \sqrt{k_2^2-\lambda^2}}{k_2} \right)^p } 
	{e^{-\sqrt{\lambda^{2}-k_2^2}(y_0-y^t_c)} e^{i\lambda (x_0-x^t_c)} } \nonumber\\
        &\times e^{(-\sqrt{k_2^2-\lambda^2} + \sqrt{k_1^2-\lambda^2} ) y_0 } { \frac{\sigma_2^t(\lambda)}{4 \pi \sqrt{\lambda^{2}-k_2^2}}  } d \lambda.
\end{align}
We set $x_0-x_c^t=2$, $y_0-y_c^t=3$, therefore the modified distance $\rho=\sqrt{2^2+3^3}$. 
The distance between the target and center of the target box is $r=1.5$. 
We neglect the propagating part in the layered medium Green's function that is very small compared with the evanescent part for large orders, 
and only consider the integral from $k_2$ to $\infty$ in the evanescent part of the basis function $\Psi_p$.
In fig.~\ref{fig:order3layer}, we show the ratio $(J_{p+1} \Psi_{p+1}) / (J_p \Psi_p) $, which clearly converges to
$r/\rho \approx 0.416$ when $p \to \infty$.

\section{Conclusion and Generalization} \label{sec:conclusion}
In this paper, we present a general numerical framework for the efficient application
of the layered media Green's function to a given density function. Instead of constructing
and compressing the matrix directly, which involves the expensive evaluations of one or more 
Sommerfeld type integrals for each matrix entry, the new algorithm considers a transformed 
matrix, so existing fast algorithms for the free-space Green's function can be readily adapted 
for better algorithm efficiency. Theoretical analysis on the convergence of the new expansions and 
alternative direction Sommerfeld integral representations to accelerate the convergence
of the numerical quadrature rules are provided and numerically validated. Similar to deriving 
the layered media Green's functions, the detailed translations, alternative direction Sommerfeld 
integral representations, and number of terms in the expansions all depend
on the geometric settings and physical parameters, especially in three dimensions. 
We have studied a few examples of such layered media Green's functions in this paper, and we are 
working on both the analysis and implementation details for other important settings from
domain applications. In particular, we are studying the {\it optimal} alternative direction
Sommerfeld integral representations and more accurate estimate of the number of expansion
terms in different scenarios. Results along these directions will be presented in the future. 

\section*{Acknowledgement}
J. Huang was supported by the NSF grant DMS1821093, and the work was finished while he was visiting professors at the 
King Abdullah University of Science and Technology, National Center for Theoretical Sciences (NCTS) in Taiwan, 
Mathematical Center for Interdisciplinary Research of Soochow University, and Institute for Mathematical Sciences 
of the National University of Singapore. M.H. Cho was supported by a grant from the Simons Foundation (No. 404499). 

\bibliography{HFMM}

\begin{thebibliography}{10}
\expandafter\ifx\csname url\endcsname\relax
  \def\url#1{\texttt{#1}}\fi
\expandafter\ifx\csname urlprefix\endcsname\relax\def\urlprefix{URL }\fi
\expandafter\ifx\csname href\endcsname\relax
  \def\href#1#2{#2} \def\path#1{#1}\fi

\bibitem{chew1995waves}
W.~C. Chew, Waves and fields in inhomogeneous media, Vol. 522, IEEE press New
  York, 1995.

\bibitem{cai2002algorithmic}
W.~Cai, Algorithmic issues for electromagnetic scattering in layered media:
  Green's functions, current basis, and fast solver, Advances in Computational
  Mathematics 16~(2-3) (2002) 157--174.

\bibitem{lai2014fast}
J.~Lai, M.~Kobayashi, L.~Greengard, A fast solver for multi-particle scattering
  in a layered medium, Optics express 22~(17) (2014) 20481--20499.

\bibitem{chen2016accurate}
D.~Chen, W.~Cai, B.~Zinser, M.~H. Cho, Accurate and efficient nystr{\"o}m
  volume integral equation method for the maxwell equations for multiple 3-d
  scatterers, Journal of Computational Physics 321 (2016) 303--320.

\bibitem{chen2018accurate}
D.~Chen, M.~H. Cho, W.~Cai, Accurate and efficient {N}ystr\"{o}m volume
  integral equation method for electromagnetic scattering of 3-{D}
  metamaterials in layered media, SIAM Journal on Scientific Computing 40~(1)
  (2018) B259--B282.

\bibitem{greengard2009fast}
L.~Greengard, D.~Gueyffier, P.-G. Martinsson, V.~Rokhlin, Fast direct solvers
  for integral equations in complex three-dimensional domains, Acta Numerica 18
  (2009) 243--275.

\bibitem{hackbusch1999sparse}
W.~Hackbusch, A sparse matrix arithmetic based on $\cal{H}$-matrices. part i:
  Introduction to $\cal{H}$-matrices, Computing 62~(2) (1999) 89--108.

\bibitem{hackbusch2000sparse}
W.~Hackbusch, B.~N. Khoromskij, A sparse $\cal{H}$-matrix arithmetic.,
  Computing 64~(1) (2000) 21--47.

\bibitem{ho2012fast}
K.~L. Ho, L.~Greengard, A fast direct solver for structured linear systems by
  recursive skeletonization, SIAM Journal on Scientific Computing 34~(5) (2012)
  A2507--A2532.

\bibitem{greengard1987fast}
L.~Greengard, V.~Rokhlin, A fast algorithm for particle simulations, Journal of
  computational physics 73~(2) (1987) 325--348.

\bibitem{cho2017efficient}
M.~H. Cho, W.~Cai, Efficient and accurate computation of electric field dyadic
  green's function in layered media, Journal of Scientific Computing 71~(3)
  (2017) 1319--1350.

\bibitem{cui1999fast}
T.~J. Cui, W.~C. Chew, Fast evaluation of sommerfeld integrals for em
  scattering and radiation by three-dimensional buried objects, IEEE
  Transactions on Geoscience and Remote Sensing 37~(2) (1999) 887--900.

\bibitem{o2014efficient}
M.~O'Neil, L.~Greengard, A.~Pataki, On the efficient representation of the
  half-space impedance green's function for the helmholtz equation, Wave Motion
  51~(1) (2014) 1--13.

\bibitem{abramowitz1966handbook}
M.~Abramowitz, I.~A. Stegun, Handbook of Mathematical Functions with Formulas,
  Graphs, and Mathematical Tables, 10th Edition, Dover, 1964.

\bibitem{rokhlin1990rapid}
V.~Rokhlin, Rapid solution of integral equations of scattering theory in two
  dimensions, Journal of Computational Physics 86~(2) (1990) 414--439.

\bibitem{crutchfield2006remarks}
W.~Crutchfield, Z.~Gimbutas, L.~Greengard, J.~Huang, V.~Rokhlin, N.~Yarvin,
  J.~Zhao, Remarks on the implementation of wideband fmm for the helmholtz
  equation in two dimensions, Contemporary Mathematics 408 (2006) 99--110.

\bibitem{greengard1998accelerating}
L.~Greengard, J.~Huang, V.~Rokhlin, S.~Wandzura, Accelerating fast multipole
  methods for the helmholtz equation at low frequencies, IEEE Computational
  Science and Engineering 5~(3) (1998) 32--38.

\bibitem{greengard1997new}
L.~Greengard, V.~Rokhlin, A new version of the fast multipole method for the
  laplace equation in three dimensions, Acta numerica 6 (1997) 229--269.

\bibitem{lai2018new}
J.~Lai, L.~Greengard, M.~O'Neil, A new hybrid integral representation for
  frequency domain scattering in layered media, Applied and Computational
  Harmonic Analysis 45~(2) (2018) 359--378.

\bibitem{cho2018heterogeneous}
M.~H. Cho, J.~Huang, D.~Chen, W.~Cai, A heterogeneous fmm for layered media
  helmholtz equation i: Two layers in r2, Journal of Computational Physics 369
  (2018) 237--251.

\bibitem{wong2001asymptotic}
R.~Wong, Asymptotic approximations of integrals, Vol.~34, SIAM, 2001.

\bibitem{cheng2006wideband}
H.~Cheng, W.~Y. Crutchfield, Z.~Gimbutas, L.~F. Greengard, J.~F. Ethridge,
  J.~Huang, V.~Rokhlin, N.~Yarvin, J.~Zhao, A wideband fast multipole method
  for the helmholtz equation in three dimensions, Journal of Computational
  Physics 216~(1) (2006) 300--325.

\bibitem{greengard2002new}
L.~F. Greengard, J.~Huang, A new version of the fast multipole method for
  screened coulomb interactions in three dimensions, Journal of Computational
  Physics 180~(2) (2002) 642--658.

\bibitem{dienstfrey2001lattice}
A.~Dienstfrey, F.~Hang, J.~Huang, Lattice sums and the two-dimensional,
  periodic green's function for the helmholtz equation, in: Proceedings of the
  Royal Society of London A: Mathematical, Physical and Engineering Sciences,
  Vol. 457, The Royal Society, 2001, pp. 67--85.

\bibitem{huang1999integral}
J.~Huang, Integral representations of harmonic lattice sums, Journal of
  Mathematical Physics 40~(10) (1999) 5240--5246.

\end{thebibliography}

\end{document}